\documentclass{ifacconf}
\usepackage{natbib}
\usepackage{amsmath,amssymb,amsfonts,bm}
\usepackage{graphicx}
\usepackage{algorithm}
\usepackage{algorithmicx,algpseudocode}
\usepackage{color}
\let\theoremstyle\relax
\usepackage{amsthm}
\usepackage[mathscr]{euscript}
\theoremstyle{plain}
\newtheorem{theorem}{Theorem}[section]
\newtheorem{lemma}[theorem]{Lemma}

\newtheorem{proposition}[theorem]{Proposition}
\theoremstyle{definition}

\theoremstyle{remark}

\begin{document}
\begin{frontmatter}

\title{Density Steering by Power Moments} 

\author[First]{Guangyu Wu} 
\author[Second]{Anders Lindquist} 

\address[First]{Department of Automation, Shanghai Jiao Tong University, Shanghai, China (e-mail: chinarustin@sjtu.edu.cn).}
\address[Second]{Department of Automation and School of Mathematical Sciences, Shanghai Jiao Tong University, Shanghai, China (e-mail: alq@kth.se)}

\begin{abstract}                
This paper considers the problem of steering an arbitrary initial probability density function to an arbitrary terminal one, where the system dynamics is governed by a first-order linear stochastic difference equation. It is a generalization of the conventional stochastic control problem where the uncertainty of the system state is usually characterized by a Gaussian distribution. We propose to use the power moments to turn the infinite-dimensional problem into a finite-dimensional one and to present an empirical control scheme. By the designed control law, the moment sequence of the controls at each time step is positive, which ensures the existence of the control for the moment system. We then realize the control at each time step as a function in analytic form by a convex optimization scheme, for which the existence and uniqueness of the solution have been proved in our previous paper \citep{wu2022non}. Two numerical examples are given to validate our proposed algorithm.
\end{abstract}

\begin{keyword}
Stochastic control; realization theory; parametric optimization; density parametrization; power moments.
\end{keyword}

\end{frontmatter}

\section{Introduction}
We are interested in the problem of steering the density function of the state where the system dynamics is governed by a discrete-time stable first-order linear stochastic difference equation. We consider the linear dynamics of the system as
\begin{equation}
x(k+1) = a(k)x(k) + u(k).
\label{uniequation}
\end{equation}
Since the system is stable and we assume $a(k)$ to be positive, we have $a(k) \in (0, 1)$. The control input to the system is defined as $u(k)$, and $x(k)$ is its state. 

Given an initial random variable $x(0)$, the density steering problems amounts to choosing a sequence of random variables $\left(u(0), u(1), \cdots, u(K-1)\right)$, so that the probability density $q_{0}$ of $x(0)$ is transferred to the density $q_{K}$ of $x(K)$ at some future time $K$.

The density steering problem has a long history and is still a cutting-edge research topic due to its significant theoretical and practical merits. For the continuous-time linear systems, Chen, Georgiou and Pavon have proposed fundamental results using the Schr\"odinger Bridge strategy for Gaussian distributions \citep{chen2015optimal, chen2015optimal2} and more general distributions \citep{chen2018optimal}. The results are extended to nonlinear continuous-time systems and hard state constraints in \cite{caluya2021reflected}. Moreover, a robust optimal density control of robotic swarms is proposed in \cite{sinigaglia2022robust}. 

For the discrete-time linear systems, the state uncertainty has usually been assumed to be Gaussian \citep{okamoto2018optimal, okamoto2019optimal, balci2020covariance}. A recent pioneering result was proposed to treat the non-Gaussian distribution steering problem by characteristic functions \citep{sivaramakrishnan2022distribution}. However, the distribution steering controller designed in the paper has an affine state feedback structure. It is then not possible for the controller to alter the function class of the distribution in the control process, which makes the controller not quite applicable to situations where the initial and terminal distributions are from different function classes, or the distributions have multiple modes. To the best of our knowledge, there has not been a complete result for the density steering problem considering discrete-time linear systems where the initial and terminal densities are non-Gaussian.

However, by generalizing the mean and covariance to all the power moments, we will have a more conceptual view of this problem. Controlling the system state as a probability density function, if only assumed to be Lebesgue integrable, is an uncountably infinite-dimensional problem, which is generally not tractable. By probability theory, we note that a distribution function can be uniquely determined by its full power moment sequence \citep{Shiryaev2016}. By controlling the full power moment sequence instead of the distribution of system state, the problem is reduced to a countably infinite-dimensional one. Then by properly truncating the first several terms of the power moment sequence for characterizing the density of the system state \citep{byrnes2006generalized, wu2022non}, the problem is now steering a truncated power moment sequence to another, which is finite-dimensional and tractable.

In this paper we provide what can be regarded as the first \textbf{computable} and \textbf{implementable} solution to the density steering problem for a discrete-time linear stochastic system in limited steps, where the specified initial and terminal density functions are \textbf{arbitrary} (only required to have first several power moments). The paper is structured as follows. In Section 2, we propose a moment system representation as a counterpart of the discrete-time linear system. There follows a formulation of the density steering by moments based on the moment system. Different from the conventional control problems, the Hankel matrices of the moments of control inputs and system states need to be positive definite, which makes it difficult to treat the control problem by prevailing methods such as optimal control. We propose an empirical scheme to treat this problem. Then we use a density parametrization algorithm proposed in our previous paper \citep{wu2022non} to realize the control inputs as analytic functions by the power moments obtained from the proposed control scheme. Last but not the least, two numerical examples show smooth transitions of the system states by the proposed algorithm. 

\section{A moment formulation of the primal problem}
In this section we treat the density steering problem formulated in Section 1. However it is not always possible to obtain a closed-form solution to this problem. If the densities are not assumed to fall within certain specific classes, the problem is intrinsically infinite-dimensional. Define the density function of the control $u(k)$ as $p_{k}(x)$, and the joint probability density function of $x(k)$ and $u(k)$ as $f_{k}\left( x, u \right)$. Then we note that the density function of $x(k+1)$ can be written as

\begin{equation}
\begin{aligned}
    q_{k+1}(t) & = \int_{\mathbb{R}} f_{k}\left( \xi, t - \xi \right)d\xi\\
    & \underline{\underline{x \perp\!\!\!\perp u}} \int_{\mathbb{R}} q_{k}\left(\frac{\xi}{a(k)}\right) p_{k}\left(t-\xi\right) d\xi\\
    & \underline{\underline{x \perp\!\!\!\perp u}} \left(q_k\left(\frac{t}{a(k)}\right) * p_k(t)\right)(t).
\end{aligned}
\label{qk1}
\end{equation}

where $x \perp\!\!\!\perp u$ denotes that $x(k)$ and $u(k)$ are independent. To tackle the density steering problem, we need to obtain a solution in analytic form of  $q_{k+1}(t)$ in \eqref{qk1}. However, except for limited classes of functions such as Gaussian distributions and trigonometric functions, this is infeasible in general. This is the main reason why in  previous results which have a similar problem setting, the examples have almost always Gaussian or trigonometric densities. This severely limits the use of these results in real applications.

There is a similar problem in non-Gaussian Bayesian filtering. In our previous results \citep{wu2022non}, we proposed a method of using the power moments to treat this intractable problem, mainly for characterizing the macroscopic property of the distributions. However, even it is theoretically feasible to characterize the probability density function by the full power moment sequence, the problem is still infinite dimensional. A common aproach is to truncate the first $2n$ moment terms \citep{byrnes2003convex, georgiou2003kullback}, which turns the problem we treat to a truncated moment problem.

By the system equation \eqref{uniequation}, the power moments of the states up to order $2n$ are written as
$$
\mathbb{E}\left[ x^{l}(k+1) \right] = \sum_{j=0}^{l}\binom{l}{j}a^{j}(k)\mathbb{E}\left[ x^{j}(k)u^{l-j}(k) \right].
$$

We note that it is difficult to treat $\mathbb{E}\left[ x^{j}(k)u^{l-j}(k) \right]$. However, we note that if $x(k)$ and $u(k)$ are independent, i.e., $\mathbb{E}\left[ x^{j}(k)u^{l-j}(k) \right] = \mathbb{E}\left[ x^{j}(k) \right]\mathbb{E}\left[ u^{i-j}(k) \right]$, the dynamics of the moments can be written as a linear matrix equation

\begin{equation}
    \mathscr{X}(k+1) = \mathscr{A}(\mathscr{U}(k))\mathscr{X}(k)+\mathscr{U}(k)
\label{momentsystem}
\end{equation}

where the state vector is composed of the power moment terms up to order $2n$, i.e.,
\begin{equation}
\mathscr{X}(k) = \begin{bmatrix}
\mathbb{E}[x(k)] & \mathbb{E}[x^{2}(k)] & \cdots & \mathbb{E}[x^{2n}(k)]
\end{bmatrix}^{T},
\label{XK}
\end{equation}
and the input vector is written as
\begin{equation}
\mathscr{U}(k) = \begin{bmatrix}
\mathbb{E}[u(k)] & \mathbb{E}[u^{2}(k)] & \cdots & \mathbb{E}[u^{2n}(k)] 
\end{bmatrix}^{T} .
\label{UK}
\end{equation}
Here 
\begin{equation}
\mathbb{E}\left[ x^{l}(k) \right] = \int_{\mathbb{R}}x^{l}q_{k}(x)dx
\label{xlK}
\end{equation}
and
$$
\mathbb{E}\left[ x^{j}(k)u^{l-j}(k) \right] = \int_{\mathbb{R}}x^{j}q_{k}(x)dx \int_{\mathbb{R}}u^{l-j}p_{k}(u)du.
$$
for $l \in \mathbb{N}_{0}$ ($\mathbb{N}_{0}$ denotes the set of all nonnegative integers), $l \leq 2n$.
Similarly we have
\begin{equation}
\mathbb{E}\left[ u^{l}(k) \right] = \int_{\mathbb{R}}u^{l}p_{k}(u)du .
\label{ulK}
\end{equation}
The matrix $\mathscr{A}(\mathscr{U}(k))$ in the system \eqref{momentsystem} can then be written as \eqref{longeq0}.

\begin{figure*}[t]
\begin{equation}
\mathscr{A}(\mathscr{U}(k))
= \begin{bmatrix}
a(k) & 0 & 0 & \cdots & 0\\ 
2a(k)\mathbb{E}[u(k)] & a^{2}(k) & 0 & \cdots & 0\\ 
3a(k)\mathbb{E}[u^{2}(k)] & 3a^{2}(k)\mathbb{E}[u(k)] & a^{3}(k) & \cdots & 0\\ 
\vdots & \vdots & \vdots & \ddots\\ 
\binom{2n}{1}a(k)\mathbb{E}[u^{2n-1}(k)] & \binom{2n}{2}a^{2}(k)\mathbb{E}[u^{2n-2}(k)] & \binom{2n}{3}a^{3}(k)\mathbb{E}[u^{2n-3}(k)] &  & a^{2n}(k)
\end{bmatrix}.
\label{longeq0}
\end{equation}
\hrulefill
\vspace*{4pt}
\end{figure*}

By using the truncated power moments to characterize the dynamics of system \eqref{uniequation} where $x(k)$ and $u(k)$ are random variables, we shall reformulate the control problem as steering the power moments of the $x(k)$ and $u(k)$. System \eqref{momentsystem} is called the moment system corresponding to system \eqref{uniequation}. The power moment steering problem is then formulated as follows.

\textit{The dynamics of the moment system is
$$
    \mathscr{X}(k+1) = \mathscr{A}(\mathscr{U}(k))\mathscr{X}(k)+\mathscr{U}(k)
$$
where $\mathscr{X}(k), \mathscr{U}(k)$ are obtained by \eqref{xlK} and \eqref{ulK}. Given an \textbf{arbitrary} initial density $q_{0}(x)$ and terminal power moments $\{\sigma_{i}\}_{i = 1:2n}$, determine the control sequence $$\left( u(0), \cdots, u(K-1)\right)$$ so that the first $2n$ order power moments of the terminal density are identical to those specified, i.e.,
\begin{equation}
\mathscr{X}(K) = \int_{\mathbb{R}} x^{l} q_{K}(x) dx = \sigma_{l}
\label{ExTl}
\end{equation}
for $l = 1, \cdots, 2n$.}

However for the moment system to control, there remains to design control laws which satisfy 

\begin{equation}
\mathbb{E}\left[ x^{j}(k)u^{l-j}(k) \right] = \mathbb{E}\left[ x^{j}(k) \right]\mathbb{E}\left[ u^{l-j}(k) \right].
\label{Expeq}
\end{equation}

To satisfy \eqref{Expeq}, a feasible control law needs to have the property that the control vector is independent of the current state vector. In the conventional feedback control law, this is hardly ever the case since the control inputs are always functions of the state vectors. However, for our density steering problem, we note that it is possible to satisfy \eqref{Expeq}, since the control inputs of the primal system, as well as the system states, are random variables. For a given system state, by drawing an i.i.d. sample from the density function of the control input, we are able to obtain a control input which is independent of the current system state. By doing this, $x(k)$ and $u(k)$ are independent, i.e., \eqref{Expeq} is satisfied. In the next section, we first propose an algorithm for steering power moments to desired ones, considering the moment system \eqref{momentsystem}.

\section{An empirical control scheme for the moment system}

In the previous section, the density steering problem was reduced to controlling the moment system corresponding to it. Then the task is now to construct an algorithm to determine a sequence of $\left( \mathscr{U}(0), \cdots, \mathscr{U}(K-1) \right)$. However, there are two main differences from the conventional control problems. First, the system matrix of the moment system is a function of the control vector in this problem. Second, the sequence of the elements of the control vector $\mathscr{U}(k)$ needs to satisfy the condition that the corresponding Hankel matrix
$$
\left[\mathscr{U}(k)\right]_{H} = \begin{bmatrix}
\mathbb{E}[u^{0}(k)] & \mathbb{E}[u^{1}(k)] & \cdots & \mathbb{E}[u^{n}(k)]\\ 
\mathbb{E}[u^{1}(k)] & \mathbb{E}[u^{2}(k)] &  & \mathbb{E}[u^{n+1}(k)]\\ 
\vdots & \vdots & \ddots & \\ 
\mathbb{E}[u^{n}(k)] & \mathbb{E}[u^{n+1}(k)] &  & \mathbb{E}[u^{2n}(k)]
\end{bmatrix}
$$

is positive definite. Here $\left[\mathscr{U}\right]_{H}$ denotes the Hankel matrix of the vector $\mathscr{U}$. We define the subspace  $\mathbb{V}^{2n}_{++} := \{ \mathscr{U} \in \mathbb{R}^{2n} \mid \left[\mathscr{U}\right]_{H} \succ 0 \}$ of $\mathbb{R}^{2n}$.

In previous results, optimal control is always used in the density steering problems. However, it doesn't work for this problem. The reason is that we always have to ensure that $\mathscr{X}(k), \mathscr{U}(k) \in \mathbb{V}^{2n}_{++}$. To the best of our knowledge, there has not been a result which is able to treat the optimal control problem constraining the states and control inputs to fall within $\mathbb{V}^{2n}_{++}$, i.e., the corresponding Hankel matrices to be positive definite. 

Let $\mathcal{U}$ be the feasible set of control sequences $\mathscr{U} := \left(\mathscr{U}(0), \cdots, \mathscr{U}(K-1) \right)$, which satisfies
$$
\sum_{k=0}^{K-1}\mathbb{E}\left[\mathscr{U}^{T}(k)\mathscr{U}(k)\right] < \infty
$$
and brings the terminal system state $x(K)$ to be distributed satisfying \eqref{ExTl}. The family $\mathcal{U}$ represents admissible control inputs which achieve the desired moment transfer. In this part of section, our goal is to obtain a proper $\mathscr{U} \in \mathcal{U}$. However in this problem, we use higer-order moments to characterize the density function of $u(k)$. Hence to determine the control inputs of the moment system by minimizing the energy effort, i.e., the second order moments, is no longer suitable. In this paper, we propose to choose the $\mathscr{U}$ by the smoothness of the transition from $\mathscr{U}(0)$ to $\mathscr{U}(K-1)$.

Even if it is very difficult to handle both the constraints on $\mathscr{U}(k)$ and $\mathscr{X}(k)$ so they fall within the set $\mathbb{V}^{2n}_{++}$ simultaneously, we note that a sub-optimal solution to the control problem can be obtained by first determining the trajectory of the state and then to obtain the control inputs corresponding to this trajectory. We first determine the trajectory of the state.

\begin{lemma}
Denote the error of the moments from the specified terminal ones as
\begin{equation}
    e(k) = \mathscr{X}(K) - \mathscr{X}(k).
\label{ek}
\end{equation} Given
$$
    e(k_{0}) = \mathscr{X}(K) - \mathscr{X}(k_{0}) \in \mathbb{V}^{2n}_{++},
$$
we have
\begin{equation}
\mathscr{X}(k) = \mathscr{X}(k_{0}) + \sum_{i = k_{0}}^{k - 1} \omega_{i} e(k_{0}) \in \mathbb{V}^{2n}_{++}
\label{Xksum}
\end{equation}
for $k = k_{0}+1, \cdots, K$ where $\omega_{i} \in \mathbb{R}_{+}$ for $i = k_{0}, \cdots, K-1$ and $\sum_{i = k_{0}}^{K - 1} \omega_{i} = 1$. Here the elements of $\mathscr{X}(K)$ are the power moments corresponding to the specified terminal density function $q_{K}(x)$.
\label{Lemma21}
\end{lemma}

\begin{proof}
The proof is straightforward. Since $\mathscr{X}(k_{0}), e(k_{0}) \in \mathbb{V}^{2n}_{++}$, we have $\left[\mathscr{X}(k_{0})\right]_{H} \succ 0$, $\left[e(k_{0})\right]_{H} \succ 0$. We note that the sum of positive definite matrices is still positive definite. Since $\omega_{i} > 0$, we have $\omega_{i} e(k_{0}) \in \mathbb{V}^{2n}_{++}$. Then $\mathscr{X}(k) \in \mathbb{V}^{2n}_{++}$.
\end{proof}

Now it remains to prove that there exists a time step $k_{0}$ at which $\mathscr{X}(K) - \mathscr{X}(k_{0}) \in \mathbb{V}^{2n}_{++}$.

\medskip

\begin{proposition}
There exists a time step $k_{0}$ which satisfies $\mathscr{X}(K) - \mathscr{X}(k_{0}) \in \mathbb{V}^{2n}_{++}$, assuming that $\mathscr{X}(k), 0 \leq k \leq k_{0}$ are uncontrolled moment states, i.e., $u(k) = 0, 0 \leq k \leq k_{0}$.
\label{proposition22}
\end{proposition}

\begin{figure*}[t]
\begin{equation}
\left[\mathscr{X}(K)- \mathscr{X}(k)\right]_{H} \\
= \begin{bmatrix}
1 & \mathbb{E}[x(K)] - \mathbb{E}[x(k)] & \cdots & \mathbb{E}[x^{n}(K)] - \mathbb{E}[x^{n}(k)]\\ 
\mathbb{E}[x(K)] - \mathbb{E}[x(k)] &  \mathbb{E}[x^{2}(K)] - \mathbb{E}[x^{2}(k)]& \cdots & \mathbb{E}[x^{n+1}(K)] - \mathbb{E}[x^{n+1}(k)]\\ 
\vdots & \vdots & \ddots & \\ 
\mathbb{E}[x^{n}(K)] - \mathbb{E}[x^{n}{k}] & \mathbb{E}[x^{n+1}(K)] - \mathbb{E}[x^{n+1}(k)] &  & \mathbb{E}[x^{2n}(K)] - \mathbb{E}[x^{2n}(k)]
\end{bmatrix}
\label{longeq1}
\end{equation}
\end{figure*}

\begin{figure*}[t]
\begin{equation}
\begin{aligned}
& \left[\mathscr{X}(K)- \mathscr{X}(k_{0})\right]_{H} = \\
&\begin{bmatrix}
1 & \mathbb{E}[x(K)] - \bar{a}\mathbb{E}[x(0)] & \cdots & \mathbb{E}[x^{n}(K)] - \bar{a}^{n}\mathbb{E}[x^{n}(0)]\\ 
\mathbb{E}[x(K)] - \bar{a}\mathbb{E}[x(0)] &  \mathbb{E}[x^{2}(K)] - \bar{a}^{2}\mathbb{E}[x^{2}(0)]& \cdots & \mathbb{E}[x^{n+1}(K)] - \bar{a}^{n+1}\mathbb{E}[x^{n+1}(0)]\\ 
\vdots & \vdots & \ddots & \\ 
\mathbb{E}[x^{n}(K)] - \bar{a}^{n}\mathbb{E}[x^{n}(0)] & \mathbb{E}[x^{n+1}(K)] - \bar{a}^{n+1}\mathbb{E}[x^{n+1}(0)] &  & \mathbb{E}[x^{2n}(K)] - \bar{a}^{2n}\mathbb{E}[x^{2n}(0)]
\end{bmatrix}
\end{aligned}
\label{longeq2}
\end{equation}
\hrulefill
\vspace*{4pt}
\end{figure*}

\begin{proof}
We write the Hankel matrix form of $\mathscr{X}(K) - \mathscr{X}(k)$ as \eqref{longeq1}. Since $u(k) = 0, 0 \leq k \leq k_{0}$, we obtain \eqref{longeq2}, where we define
$$
\bar{a} = \prod_{i = 0}^{k - 1}a(i).
$$
Now it remains to prove that $\exists k_{0}$, $\left[\mathscr{X}(K)- \mathscr{X}(k_{0})\right]_{H} \succ 0$ with $u(k) = 0, 0 \leq k \leq k_{0}$. By definition of  positive definiteness, this is equivalent to proving that each leading principal minor, the determinant of leading principal submatrix, is positive.

Here we denote the $i_\text{th}$-order leading principal submatrix of $\left[\mathscr{X}(K)- \mathscr{X}(k)\right]_{H}$ as $H_{i}(k)$, and the corresponding minor as $\operatorname{det}(H_{i}(k))$. We note that each $\operatorname{det}(H_{i}(k))$ is a polynomial of $\bar{a} = \prod_{i = 0}^{k - 1}a(i)$, of which the degree is even. Therefore, if there exists no real zero for all the $\operatorname{det}(H_{i}(k))$, all $k_{0} \in \mathbb{N}_{0}$ satisfies $\left[\mathscr{X}(K)- \mathscr{X}(k_{0})\right]_{H} \succ 0$. Now we consider the case that $\operatorname{det}(H_{i}(k))$ has at least a real zero in $(0, 1)$. We note that $\operatorname{det}(H_{i}(k_{0})) > 0$ with $k_{0} \rightarrow +\infty$. Let $\breve{k}_{i}$ be the smallest integer greater than the largest zero of the polynomial $\operatorname{det}(H_{i}(k))$. By the continuity of $\operatorname{det}(H_{i}(k))$, we have that $\operatorname{det}(H_{i}(k)) > 0, k \in ( \breve{k}_{i}, +\infty)$.

Let $
\breve{k} = \max_{i}(\breve{k}_{i})$.
With $k_{0} > \breve{k}$, we have $H_{i}(k_{0}) \succ 0$, for $i = 1, \cdots, n$. Therefore we have $\left[\mathscr{X}(K)- \mathscr{X}(k_{0})\right]_{H} \succ 0$, which ensures the positiveness of all $H_{i}(k_{0})$ and completes the proof.
\end{proof}

By Proposition \ref{proposition22}, it is possible to choose a time step $k_{0}$ which satisfies $\mathscr{X}(K) - \mathscr{X}(k_{0}) \in \mathbb{V}^{2n}_{++}$. We assume that the system is uncontrolled before $k_{0}$, i.e. $u(k)=0, k \leq k_{0}$. From step $k_{0}$, we impose controls on the system. Lemma \ref{Lemma21} has proved the positiveness of $\mathscr{X}(k), k = k_{0}, \cdots, K$. Therefore it remains to determine the parameters $\omega_{k}, k = k_{0}, \cdots, K - 1$ and the corresponding control inputs $\mathscr{U}(k)$.

This is a non-trivial problem. We give an empirical scheme to treat it. To obtain a relatively smooth transition of states, it is desired that $\omega_{i}$ for $i = k_{0}, \cdots, K-1$ are close to each other. We are usually able to choose
$$
\omega_{k_{0}} = \cdots = \omega_{K-1} = \frac{1}{K-k_{0}}
$$
After that the parameters $\omega_{i}$ are determined. Then the control inputs of the moment system $\mathscr{U}(i)$ for $i = k_{0}, \cdots, K-1$ can be calculated by solving the equation \eqref{momentsystem}, provided with $\mathscr{X}(k), k = k_{0}+1, \cdots, K$ calculated by \eqref{Xksum}. The obtained moment sequence $\left(  \mathscr{U}(0), \cdots, \mathscr{U}(K) \right) \in \mathcal{U}$ is then a solution to the moment steering problem.

However sometimes the control inputs $\mathscr{U}(k) \notin \mathbb{V}^{2n}_{++}$ by choosing the $\omega_{i}$ for $i = k_{0}, \cdots, K-1$ to be all equal. It usually happens when the specified initial or terminal density has several modes (peaks). If so, we can choose a larger $\omega_{0}$ or $\omega_{K-1}$. 

In conclusion, we have proposed to use the moment system to steer the probability density functions. And an empirical control law has been proposed which ensures the existence of $u(k)$. However by the proposed control scheme for the moment system, we are only able to obtain the power moments of the control inputs $u(k)$. We need to obtain the $u(k)$ given its power moments, which we call realization of the control inputs.

\section{Realization of the control inputs}
In this section we shall realize the probability density of $u(k)$ given the power moments of the designed controls $\mathscr{U}(k)$ for the moment system. 

For the sake of simplicity, we omit $k$ if there is no ambiguity in the following part of this section. The problem now becomes that of proposing an algorithm which  estimates the probability density for which the power moments are as specified.

A convex optimization scheme for density estimation by the Kullback-Leibler distance has been proposed in \citep{wu2022non} considering the Hamburger moment problem. We adopt this strategy in this paper for treating the realization of the control inputs. The procedure is as follows.

Let $\mathcal{P}$ be the space of probability density functions on the real line with support there, and let $\mathcal{P}_{2n}$ be the subset of all $p \in \mathcal{P}$ which have at least $2n$ finite moments (in addition to $\mathbb{E}[u^{0}(k)]$, which of course is 1). The Kullback-Leibler distance is then defined as
\begin{equation}
\mathbb{K} \mathbb{L}(r \| p)=\int_{\mathbb{R}} r(u) \log \frac{r(u)}{p(u)} d u
\label{Hellinger}
\end{equation} 
where $r$ is an arbitrary probability density in $\mathcal{P}$. We define the linear integral operator $\Gamma$ as
$$
\Gamma: p(u) \mapsto \Sigma=\int_{\mathbb{R}} G(u) p(u) G^{T}(u) d u,
$$
where $p(u)$ belongs to the space $\mathcal{P}_{2n}$. Here
$$
G(u)= \begin{bmatrix}
1 & u & \cdots & u^{n-1} & u^{n}
\end{bmatrix}^{T}
$$
and
$$
\Sigma = \begin{bmatrix}
1 & \mathbb{E}[u] & \cdots & \mathbb{E}[u^{n}]\\ 
\mathbb{E}[u] &  \mathbb{E}[u^{2}]& \cdots & \mathbb{E}[u^{n+1}]\\ 
\vdots & \vdots & \ddots & \\ 
\mathbb{E}[u^{n}] & \mathbb{E}[u^{n+1}] &  & \mathbb{E}[u^{2n}(K)]
\end{bmatrix}
$$
where $\mathbb{E}[u^{i}], i = 1, \cdots, 2n$ are the elements of the designed control $\mathscr{U}$. Moreover, since $\mathcal{P}_{2n}$ is convex, then so is $\operatorname{range}(\Gamma)=\Gamma\mathcal{P}_{2n}$.

We let
$$
\mathcal{L}_{+}:=\left\{\Lambda \in \operatorname{range}(\Gamma) \mid G(u)^{T} \Lambda G(u)>0, x \in \mathbb{R}\right\}.
$$
Given any $r \in \mathcal{P}$ and any $\Sigma \succ 0$, there is a unique $\hat{p} \in \mathcal{P}_{2 n}$ that minimizes \eqref{Hellinger} subject to $\Gamma(\hat{p})=\Sigma$, namely
\begin{equation}
\hat{p}=\frac{r}{G^{T} \hat{\Lambda} G} 
\label{hatp}
\end{equation}
where $\hat{\Lambda}$ is the unique solution to the problem of minimizing
\begin{equation}
\mathbb{J}_{r}(\Lambda):=\operatorname{tr}(\Lambda \Sigma)-\int_{\mathbb{R}} r(u) \log \left[G(u)^{T} \Lambda G(u)\right] d u
\label{Jr}
\end{equation}
over all $\Lambda \in \mathcal{L}_{+}$   \citep{wu2022non}. Here $\operatorname{tr}(M)$ denotes the trace of the matrix $M$.

Then the density estimation is formulated as a convex optimization problem. The map $\Lambda \mapsto \Sigma$ has been proved to be homeomorphic, which ensures the existence and uniqueness of the solution to the realization of control inputs \citep{wu2022non}. Unlike other moment methods, the power moments of our proposed density estimate are exactly identical to those specified, which makes it a satisfactory approach for realization of the control inputs \citep{wu2022non}. Since the prior density $r(u)$ and the density estimate $\hat{p}(u)$ are both supported on $\mathbb{R}$, $r(u)$ can be chosen as a Gaussian distribution (or a Cauchy distribution if $\hat{p}(u)$ is assumed to be heavy-tailed).

\section{Numerical examples}
In this section, we perform numerical simulations on mixtures of density functions supported on $\mathbb{R}$ to validate our proposed algorithms.

\begin{figure}[htbp]
\centering
\includegraphics[scale=0.32]{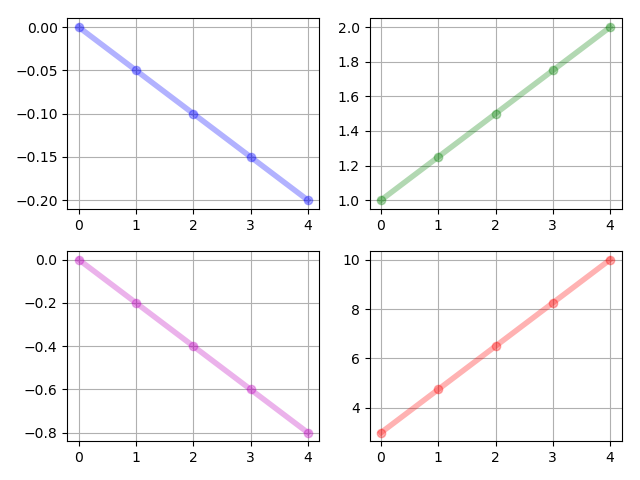}
\centering
\caption{$\mathscr{X}(k)$ at time steps $k = 0, 1, 2, 3, 4$. The upper left figure shows $\mathbb{E}\left[ x(k)\right]$. The upper right one shows $\mathbb{E}\left[ x^{2}(k)\right]$. The lower left one shows $\mathbb{E}\left[ x^{3}(k)\right]$ and the lower right one shows $\mathbb{E}\left[ x^{4}(k)\right]$.}
\label{fig1}
\end{figure}

\begin{figure}[htbp]
\centering
\includegraphics[scale=0.32]{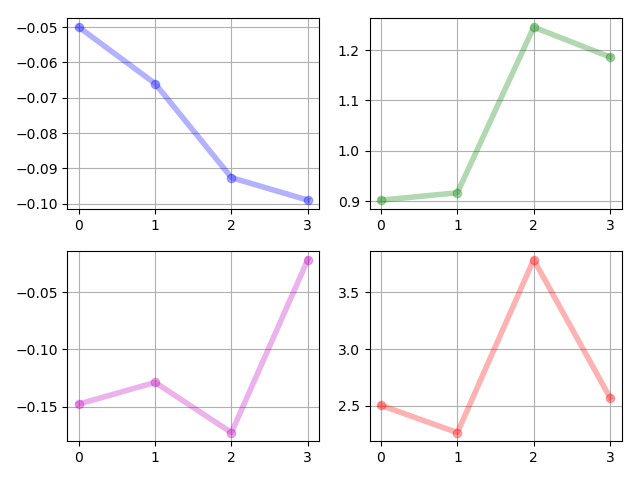}
\centering
\caption{$\mathscr{U}(k)$ at time steps $k = 0, 1, 2, 3$. The upper left figure shows $\mathbb{E}\left[ u(k)\right]$. The upper right one shows $\mathbb{E}\left[ u^{2}(k)\right]$. The lower left one shows $\mathbb{E}\left[ u^{3}(k)\right]$ and the lower right one shows $\mathbb{E}\left[ u^{4}(k)\right]$.}
\label{fig2}
\end{figure}

\begin{figure}[htbp]
\centering
\includegraphics[scale=0.32]{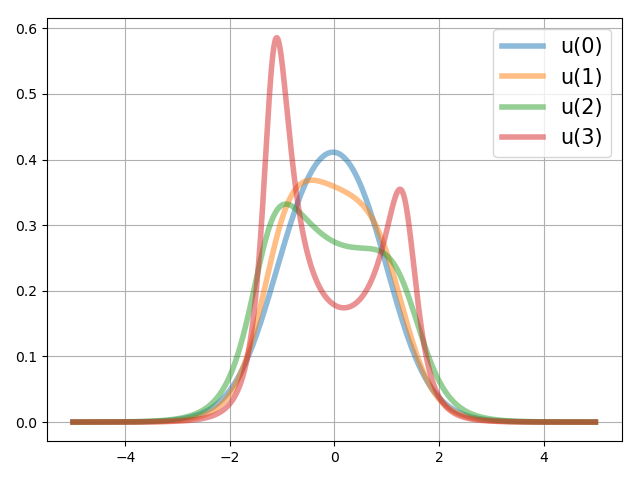}
\centering
\caption{Probability densities of the realized control inputs $u(k)$ by $\mathscr{U}(k)$ for $k = 0, 1, 2, 3$, which are obtained by our proposed control scheme.}
\label{fig3}
\end{figure}

In Example 1, we simulate a problem which is to steer a Gaussian density to a mixture of Gaussians in four steps. The initial Gaussian density is chosen as 
\begin{equation}
     q_{0}(x) = \frac{1}{\sqrt{2\pi}}e^{\frac{x ^{2}}{2}},
\label{q01}
\end{equation}
and the terminal density is specified as
\begin{equation}
     q_{K}(x) = \frac{0.4}{\sqrt{2\pi}}e^{\frac{(x - 1) ^{2}}{2}} + \frac{0.6}{\sqrt{2\pi}}e^{\frac{(x + 1) ^{2}}{2}}.
\label{qt1}
\end{equation}
The system parameters $a(k), k = 0, \cdots, 3$ are i.i.d. samples drawn from the uniform distribution $U[0.5, 0.7]$. In this example, the dimension of each $\mathscr{U}(k)$ is $4$, i.e., power moments up to order $4$ are used for realizing the control inputs $u(k)$ for $k = 0, \cdots, 3$. The states of the moment system, i.e., $\mathscr{X}(k)$ for $k = 0, 1, 2, 3, 4$ are given in Figure \ref{fig1}. The controls of the moment system, i.e., $\mathscr{U}(k)$ for $k = 0, 1, 2, 3$ are given in Figure \ref{fig2}. We note that by our proposed algorithm, $\mathscr{X}(k), \mathscr{U}(k) \in \mathbb{V}^{2n}_{++}$, which makes it possible for us to realize the controls. The realized control inputs are given in Figure \ref{fig3}. The transition of the control inputs is smooth, which is a satisfactory result. 

\begin{figure}[htbp]
\centering
\includegraphics[scale=0.32]{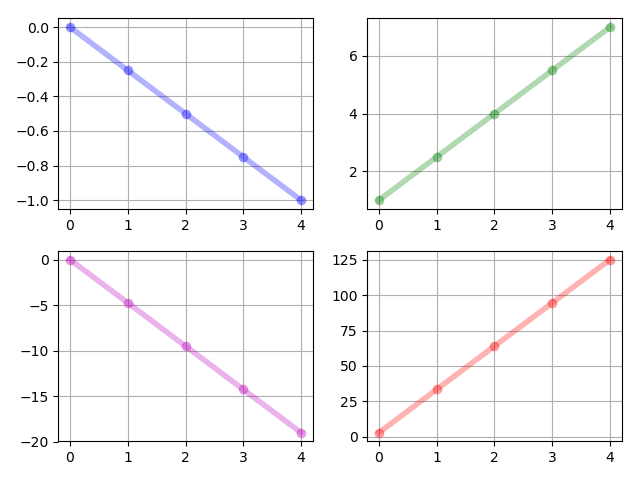}
\centering
\caption{$\mathscr{X}(k)$ at time steps $k = 0, 1, 2, 3, 4$.}
\label{fig4}
\end{figure}

\begin{figure}[htbp]
\centering
\includegraphics[scale=0.32]{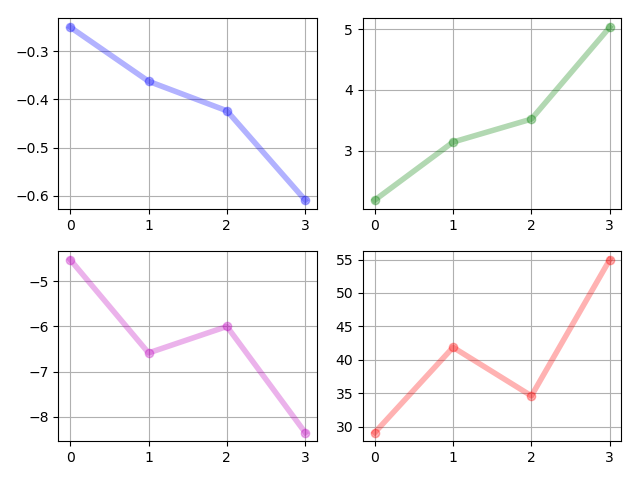}
\centering
\caption{$\mathscr{U}(k)$ at time steps $k = 0, 1, 2, 3$.}
\label{fig5}
\end{figure}

\begin{figure}[htbp]
\centering
\includegraphics[scale=0.32]{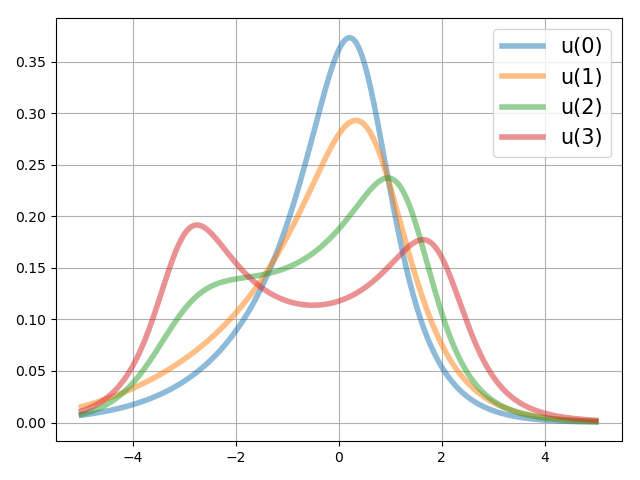}
\centering
\caption{Probability densities of the realized control inputs $u(k)$ for $k = 0, 1, 2, 3$.}
\label{fig6}
\end{figure}

In Example 2, we simulate a steering problem in four steps where the initial density function is a Gaussian and the terminal density function is a mixture of Laplacians with two modes. The initial one is chosen as
$$
     q_{0}(x) = \frac{1}{\sqrt{2\pi}}e^{\frac{x ^{2}}{2}}.
$$
and the terminal one is specified as
\begin{equation}
     q_{K}(x) = \frac{0.5}{2}e^{\left| x - 1 \right|} + \frac{0.5}{2}e^{-\left| x + 3\right|}.
\label{qt2}
\end{equation}
The system parameters $a(k), k = 0, \cdots, 3$ are i.i.d. samples drawn from the uniform distribution $U[0.5, 0.7]$. The dimension of each $\mathscr{U}(k)$ is $4$. The states of the moment system, i.e., $\mathscr{X}(k)$ for $k = 0, 1, 2, 3, 4$ are given in Figure \ref{fig4}. The controls of the moment system, i.e., $\mathscr{U}(k)$ for $k = 0, 1, 2, 3$ are given in Figure \ref{fig5}. The realized controls in Figure \ref{fig6} also show that the transition of the control inputs is smooth, even the specified terminal density has two modes, which are Laplacians.

\section{Conclusions}

In this paper, we treat the problem of steering an arbitrary probability density function to another one for the discrete-time stable first-order stochastic linear system, for the first time in the literature. The infinite-dimensional primal problem is turned to a finite-dimensional one by formulating it as a moment-steering problem. Then an empirical control scheme is proposed with the constraints that all Hankel matrices of the states and the controls are all positive definite. A realization of the controls is also given. Two numerical examples validate our algorithm where the transitions of the system states are smooth.

In future work, we are interested in extending the results of this paper to multi-dimensional systems. Moreover, we would like to discover the relationship between the moment method and the Schr\"odinger Bridge, since the state transitions for both of them are smooth.

\bibliography{ifacconf}             
\end{document}